\documentclass[a4paper,11pt]{amsart}
\usepackage{graphicx}
\usepackage{mathptmx}
\usepackage{amsmath}
\usepackage{amssymb}
\usepackage{enumitem}
\usepackage{xcolor}
\usepackage{xparse}
\NewDocumentCommand{\eulerian}{omm}
 {%
  \genfrac<>{0pt}{}{#2}{#3}%
  \IfValueT{#1}{_{\!#1}}%
 }

 \newmuskip\pFqmuskip

\newcommand*\pFq[6][8]{%
  \begingroup 
  \pFqmuskip=#1mu\relax
  \mathchardef\normalcomma=\mathcode`,
  \mathcode`\,=\string"8000
  \begingroup\lccode`\~=`\,
  \lowercase{\endgroup\let~}\pFqcomma
  {}_{#2}F_{#3}{\left(\genfrac..{0pt}{}{#4}{#5}\bigg|#6\right)}%
  \endgroup
}
\newcommand{\pFqcomma}{{\normalcomma}\mskip\pFqmuskip}

\newtheorem{theorem}{Theorem}
\newtheorem{lemma}[theorem]{Lemma}

\newtheorem{proposition}[theorem]{Proposition}

\newtheorem{remark}[theorem]{Remark}

\begin{document}

\title[Generalized degenerate Bernoulli numbers and polynomials]{Generalized degenerate Bernoulli numbers and polynomials arising from Gauss hypergeometric function}

\author{Taekyun  Kim $^{1,\dagger}$}
\address{$^{1}$ Department of Mathematics, Kwangwoon University, Seoul 139-701, Republic of Korea}
\email{tkkim@kw.ac.kr$^{\dagger}$, luciasconstant@kw.ac.kr$^{\ddagger}$, gksdud213@kw.ac.kr}

\author{DAE SAN KIM $^{2}$*}
\address{$^{2}$ Department of Mathematics, Sogang University, Seoul 121-742, Republic of Korea}
\email{dskim@sogang.ac.kr*}

\author{Lee-Chae Jang $^{3}$**}
\address{$^{3}$ Graduate School of Education, Konkuk University, Seoul 143-701, Republic of Korea}
\email{Lcjang@konkuk.ac.kr**}

\author{Hyunseok Lee $^{1,\ddagger}$}

\author{Hanyoung Kim $^{1,\dagger\dagger}$}

\subjclass[2010]{11B68; 11B73; 11B83; 33C05}
\keywords{generalized degenerate Bernoulli numbers; generalized degenerate Bernoulli polynomials; degenerate type Eulerian numbers}

\maketitle

\begin{abstract}
In a previous paper, Rahmani introduced a new family of $p$-Bernoulli numbers and polynomials by means of the Gauss hypergeometric function. Motivated by this paper and as a degenerate version of those numbers and polynomials, we introduce the generalized degenerate Bernoulli numbers and polynomials again by using the Gauss hypergeometric function. In addition, we introduce the degenerate type Eulerian numbers as a degenerate version of Eulerian numbers. For the generalized degenerate Bernoulli numbers, we express them in terms of the degenerate Stirling numbers of the second kind, of the degenerate type Eulerian numbers, of the degenerate $p$-Stirling numbers of the second kind and of an integral on the unit interval. As to the generalized degenerate Bernoulli polynomials, we represent them in terms of the degenerate Stirling polynomials of the second kind.
\end{abstract}

\section{Introduction}
As the first degenerate versions of some special numbers, Carlitz introduced the degenerate Stirling, Bernoulli and Euler numbers in [3].
In recent years, degenerate versions of many special polynomials and numbers have been investigated by means of various different tools including generating functions, combinatorial methods, umbral calculus, $p$-adic analysis, differential equations, special functions, probability theory and analytic number theory. Here we would like to remark that studying degenerate versions of some special polynomials and numbers has yielded many interesting arithmetic and combinatorial results (see [7-13 and references therein]) and has potential to find many applications to diverse areas in science and engineering as well as in mathematics. For example, it was shown in [10,11] that both the degenerate $\lambda$-Stirling polynomials of the
second kind and the $r$-truncated degenerate $\lambda$-Stirling polynomials of the second kind appear in the
expressions of the probability distributions of appropriate random variables. Also, we would like to emphasize that studying degenerate versions is applied not only to polynomials but also to transcendental functions. Indeed, the degenerate gamma functions were introduced and some interesting results were derived in [9]. \par
In [14], Rahmani introduced a new family of $p$-Bernoulli numbers and polynomials by means of the Gauss hypergeometric function which reduce to the classical Bernoulli numbers and polynomials for $p=0$. Motivated by that paper and as a degenerate version of those numbers and polynomials, in this paper we introduce the generalized degenerate Bernoulli numbers and polynomials again in terms of the Gauss hypergeometric function which reduce to the Carlitz degenerate Bernoulli numbers and polynomials for $p=0$. In addition, we introduce the degenerate type Eulerian numbers as a degenerate version of Eulerian numbers. The aim of this paper is to study the generalized degenerate Bernoulli numbers and polynomials and to show their connections to other special numbers and polynomials.
Among other things, for the generalized degenerate Bernoulli numbers we express them in terms of the degenerate Stirling numbers of the second kind , of the degenerate type Eulerian numbers, of the degenerate $p$-Stirling numbers of the second kind and of an integral on the unit interval. As to the generalized degenerate Bernoulli polynomials, we represent them in terms of the degenerate Stirling polynomials of the second kind.
For the rest of this section, we recall the necessary facts that are needed throughout this paper.
\vspace{0.2cm}

For any $\lambda\in\mathbb{R}$, the degenerate exponential functions are defined by
\begin{equation}
e_{\lambda}^{x}(t)=\sum_{n=0}^{\infty}(x)_{n,\lambda}\frac{t^{n}}{n!},\quad e_{\lambda}(t)=e_{\lambda}^{1}(t),\quad(\mathrm{see}\ [6,9]),\label{1}
\end{equation}
where $(x)_{0,\lambda}=1,\ (x)_{n,\lambda}=x(x-\lambda)\cdots(x-(n-1)\lambda)$, $(n\ge 1)$.
Note that $\displaystyle\lim_{\lambda\rightarrow 0}e^{x}_{\lambda}(t)=e^{xt}\displaystyle$. \par
Let $\log_{\lambda}(t)$ be the compositional inverse function of $e_{\lambda}(t)$ with $\log_{\lambda}\big(e_{\lambda}(t)\big)=e_{\lambda}\big(\log_{\lambda}(t)\big)=t$. Then we have
\begin{equation}
\log_{\lambda}(1+t)=\sum_{n=1}^{\infty}\lambda^{n-1}(1)_{n,1/\lambda}\frac{t^{n}}{n!},\quad(\mathrm{see}\ [7]).\label{2}	
\end{equation}
In [7], the degenerate Stirling numbers of the first kind are defined by
\begin{equation}
(x)_{n}=\sum_{l=0}^{n}S_{1,\lambda}(n,l)(x)_{l,\lambda},\quad(n\ge 0),\label{3}	
\end{equation}
where $(x)_{0}=1,\ (x)_{n}=x(x-1)(x-2)\cdots(x-n+1)$, $(n\ge 1)$. \par
As the inversion formula of $\eqref{3}$, the degenerate Stirling numbers of the second kind are defined by
\begin{equation}
(x)_{n,\lambda}=\sum_{k=0}^{n}S_{2,\lambda}(n,k)(x)_{k},\quad(n\ge 0),\quad(\mathrm{see}\ [7]).\label{4}
\end{equation}
From \eqref{3} and \eqref{4}, we note that
\begin{equation}
\frac{1}{k!}\big(\log_{\lambda}(1+t)\big)^{k}=\sum_{n=k}^{\infty}S_{1,\lambda}(n,k)\frac{t^{n}}{n!},\label{5}
\end{equation}
and
\begin{equation}
\frac{1}{k!}\big(e_{\lambda}(t)-1\big)^{k}=\sum_{n=k}^{\infty}S_{2,\lambda}(n,k)\frac{t^{n}}{n!},\quad(k\ge 0),\quad(\mathrm{see}\ [7]).\label{6}
\end{equation}
It is well known that the Gauss hypergeometric function is given by
\begin{equation}
\pFq{2}{1}{a,b}{c}{x}=\sum_{k=0}^{\infty}\frac{\langle a\rangle_{k}\langle b\rangle_{k}}{\langle c\rangle_{k}}\frac{x^{k}}{k!},\quad(\mathrm{see}\ [1,2,12]),\label{7}
\end{equation}
where $\langle a\rangle_{0}=1,\ \langle a\rangle_{k}=a(a+1)\cdots(a+k-1),\ (k\ge 1)$. \par
The Pfaff's transformation formula is given by
\begin{equation}
\pFq{2}{1}{a,b}{c}{x}=(1-x)^{-a}\pFq{2}{1}{a,c-b}{c}{\frac{x}{x-1}},\quad(\mathrm{see}\ [1,2]),\label{8}
\end{equation}
and the Euler's transformation formula is given by
\begin{equation}
\pFq{2}{1}{a,b}{c}{x}=(1-x)^{c-a-b}\pFq{2}{1}{c-a,c-b}{c}{x},\quad(\mathrm{see}\ [1,2]).\label{9}
\end{equation}
The Eulerian number $\eulerian{n}{k}$ is the number of permutation $\{1,2,3,\dots,n\}$ having $k$ permutation ascents. The Eulerian numbers are given explicitly by the finite sum
\begin{equation}
\eulerian{n}{k}=\sum_{j=0}^{k+1}(-1)^{j}\binom{n+1}{j}(k-j+1)^{n},\quad(n, k \ge 0,n \ge k), \label{10}
\end{equation}
and
\begin{equation}
\sum_{k=0}^{n}\eulerian{n}{k}=n!,\quad(\mathrm{see}\ [4,5]).\label{11}
\end{equation}
For $n,m\ge 0$, we have
\begin{equation}
\eulerian{n}{m}=\sum_{k=0}^{n-m}S_{2}(n,k)\binom{n-k}{m}(-1)^{n-k-m}k!,\quad(\mathrm{see}\ [5]),\label{12}
\end{equation}
and
\begin{equation}
x^{n}=\sum_{k=0}^{n}\eulerian{n}{k}\binom{x+k}{n},\quad(\mathrm{see}\ [4,5]).\label{13}
\end{equation}
Recently, the degenerate Stirling polynomials of the second kind are defined by
\begin{equation}
\frac{1}{k!}\big(e_{\lambda}(t)-1\big)^{k}e_{\lambda}^{x}(t)=\sum_{n=k}^{\infty}S_{2,\lambda}(n,k|x)\frac{t^{n}}{n!},\quad(k\ge 0),\quad(\mathrm{see}\ [8]).\label{14}
\end{equation}
Thus, by \eqref{14}, we get
\begin{align}
S_{2,\lambda}(n,k|x)\ &=\ \sum_{l=k}^{n}\binom{n}{l}S_{2,\lambda}(l,k)(x)_{n-l,\lambda},\quad(\mathrm{see}\ [8]), \label{15} \\
&=\ \sum_{l=0}^{n}\binom{n}{l}S_{2,\lambda}(l,k)(x)_{n-l,\lambda},\quad(n\ge 0).\nonumber
\end{align}
For $x=0$, $S_{2,\lambda}(n,k)=S_{2,\lambda}(n,k|0)$, $(n,k\ge 0, n\ge k)$, are called the degenerate Stirling numbers of the second kind.\par
Carlitz introduced the degenerate Bernoulli polynomials given by
\begin{equation}
\frac{t}{e_{\lambda}(t)-1}e_{\lambda}^{x}(t)=\sum_{n=0}^{\infty}\beta_{n,\lambda}(x)\frac{t^{n}}{n!},\quad(\mathrm{see}\ [3]). \label{16}
\end{equation}
When $x=0$, $\beta_{n,\lambda}=\beta_{n,\lambda}(0)$, $(n\ge 0)$, are called the degenerate Bernolli numbers. \par

\section{Generalized degenerate Bernoulli numbers}
By \eqref{1} and \eqref{2}, we get
\begin{align}
\frac{t}{e_{\lambda}(t)-1}\ &=\ \frac{1}{e_{\lambda}(t)-1}\sum_{n=1}^{\infty}\lambda^{n-1}(1)_{n,1/\lambda}\frac{1}{n!}\big(e_{\lambda}(t)-1\big)^{n}\label{17} \\
&=\ \sum_{k=0}^{\infty}\frac{\lambda^{k}(1)_{k+1,1/\lambda}}{k+1}\cdot\frac{1}{k!}\big(e_{\lambda}(t)-1\big)^{k}\nonumber \\
&=\ \sum_{k=0}^{\infty}\frac{\lambda^{k}(1)_{k+1,1/\lambda}}{k+1}\sum_{n=k}^{\infty}S_{2,\lambda}(n,k)\frac{t^{n}}{n!} \nonumber \\
&=\ \sum_{n=0}^{\infty}\bigg(\sum_{k=0}^{n}\frac{\lambda^{k}(1)_{k+1,1/\lambda}}{k+1}S_{2,\lambda}(n,k)\bigg)\frac{t^{n}}{n!}.\nonumber
\end{align}
Therefore, by \eqref{16} and \eqref{17}, we obtain the following theorem.
\begin{theorem}
	For $n\ge 0$, we have
	\begin{displaymath}
		\beta_{n,\lambda}=\sum_{k=0}^{n}\frac{\lambda^{k}(1)_{k+1,1/\lambda}}{k+1}S_{2,\lambda}(n,k).
	\end{displaymath}
\end{theorem}
Replacing $t$ by $\log_{\lambda}(1+t)$ in \eqref{16}, we get
\begin{align}
\frac{\log_{\lambda}(1+t)}{e_{\lambda}(\log_{\lambda}(1+t))-1}\ &=\ \sum_{k=0}^{\infty}\beta_{k,\lambda}\frac{1}{k!}\big(\log_{\lambda}(1+t)\big)^{k}\label{18} \\
&=\ \sum_{k=0}^{\infty}\beta_{k,\lambda}\sum_{n=k}^{\infty}S_{1,\lambda}(n,k)\frac{t^{n}}{n!} \nonumber \\
&=\ \sum_{n=0}^{\infty}\bigg(\sum_{k=0}^{n}S_{1,\lambda}(n,k)\beta_{k,\lambda}\bigg)\frac{t^{n}}{n!}.\nonumber
\end{align}
On the other hand, by \eqref{2}, we get
\begin{align}
\frac{\log_{\lambda}(1+t)}{e_{\lambda}(\log_{\lambda}(1+t))-1}\ &=\ \frac{1}{t}\log_{\lambda}(1+t)\ =\ \frac{1}{t}\sum_{n=1}^{\infty}\lambda^{n-1}(1)_{n,1/\lambda}\frac{t^{n}}{n!} \label{19} \\
&=\ \sum_{n=0}^{\infty}\frac{\lambda^{n}(1)_{n+1,1/\lambda}}{n+1}\frac{t^{n}}{n!}.\nonumber
\end{align}
Therefore, by \eqref{18} and \eqref{19}, we obtain the following theorem.
\begin{theorem}
	For $n\ge 0$, we have
	\begin{displaymath}
		\sum_{k=0}^{n}S_{1,\lambda}(n,k)\beta_{k,\lambda}=\frac{1}{n+1}\lambda^{n}(1)_{n+1,1/\lambda}.
	\end{displaymath}
\end{theorem}
From \eqref{16} and \eqref{17}, we note that
\begin{align}
\sum_{n=0}^{\infty}\beta_{n,\lambda}\frac{t^{n}}{n!}\ &=\ \frac{1}{e_{\lambda}(t)-1}\sum_{n=1}^{\infty}\lambda^{n-1}(1)_{n,1/\lambda}\frac{1}{n!}\big(e_{\lambda}(t)-1\big)^{n}\label{20} \\
&=\ \sum_{n=0}^{\infty}\frac{(-1)^{n}(1)_{n+1,1/\lambda}\lambda^nn!}{(n+1)!}\frac{(1-e_{\lambda}(t))^{n}}{n!}\nonumber \\
&=\ \sum_{n=0}^{\infty}\frac{\langle 1-\lambda\rangle_{n}\langle 1\rangle_{n}}{\langle 2\rangle_{n}}\frac{(1-e_{\lambda}(t))^{n}}{n!}\ \nonumber\\
&=\ \pFq{2}{1}{1-\lambda,1}{2}{1-e_{\lambda}(t)}.\nonumber	
\end{align}
In view of \eqref{20}, we may consider the {\it{generalized degenerate Bernoulli numbers}} given in terms of  Gauss hypergeometric function by
\begin{equation}
\pFq{2}{1}{1-\lambda,1}{p+2}{1-e_{\lambda}(t)}=\sum_{n=0}^{\infty}\beta_{n,\lambda}^{(p)}\frac{t^{n}}{n!},\label{21}
\end{equation}
where $p\in\mathbb{Z}$ with $p\ge -1$. When $p=0$, $\beta_{n,\lambda}^{(0)}=\beta_{n,\lambda},\ (n\ge 0)$. \par
Let us take $p=-1$ in \eqref{21}. Then we have
\begin{align}
\sum_{n=0}^{\infty}\beta_{n,\lambda}^{(-1)}\frac{t^{n}}{n!}\ &=\ \pFq{2}{1}{1-\lambda,1}{1}{1-e_{\lambda}(t)}\nonumber \\
&=\ \sum_{n=0}^{\infty}\frac{(\lambda-1)_{n}}{n!}\big(e_{\lambda}(t)-1\big)^{n}\ =\ \sum_{n=0}^{\infty}\binom{\lambda-1}{n}\big(e_{\lambda}(t)-1\big)^{n}\label{22}\\
&=\ e_{\lambda}^{\lambda-1}(t)\ =\ \sum_{n=0}^{\infty}(\lambda-1)_{n,\lambda}\frac{t^{n}}{n!}.\nonumber
\end{align}
By comparing the coefficients on the both sides of \eqref{22}, we get
\begin{equation}
\beta_{n,\lambda}^{(-1)}=(\lambda-1)_{n,\lambda},\quad(n\ge 0). \label{23}	
\end{equation}
From \eqref{21}, we note that
\begin{align}
\sum_{n=0}^{\infty}\beta_{n,\lambda}^{(p)}\frac{t^{n}}{n!}\ &=\ \pFq{2}{1}{1-\lambda,1}{p+2}{1-e_{\lambda}(t)}\ =\ \sum_{k=0}^{\infty}\frac{\langle 1-\lambda\rangle_{k}\langle 1\rangle_{k}}{\langle p+2\rangle_{k}}\frac{(1-e_{\lambda}(t))^{k}}{k!}\label{24} \\
&=\ (p+1)!\sum_{k=0}^{\infty}\frac{\lambda^{k}(1)_{k+1,1/\lambda}k!}{(p+k+1)!}\frac{1}{k!}\big(e_{\lambda}(t)-1\big)^{k} \nonumber \\
&=\ (p+1)!\sum_{k=0}^{\infty}\frac{\lambda^{k}(1)_{k+1,1/\lambda}k!}{(p+k+1)!}\sum_{n=k}^{\infty}S_{2,\lambda}(n,k)\frac{t^{n}}{n!} \nonumber \\
&=\ \sum_{n=0}^{\infty}\bigg(\sum_{k=0}^{n}\frac{\lambda^{k}(1)_{k+1,1/\lambda}}{\binom{p+k+1}{p+1}}S_{2,\lambda}(n,k)\bigg)\frac{t^{n}}{n!}.\nonumber
\end{align}
Therefore, by comparing the coefficients on both sides of \eqref{24}, we obtain the following theorem.
\begin{theorem}
	For $n\ge 0$ and $	p\ge -1$, we have
	\begin{displaymath}
		\beta_{n,\lambda}^{(p)}= \sum_{k=0}^{n}\frac{\lambda^{k}(1)_{k+1,1/\lambda}}{\binom{p+k+1}{p+1}}S_{2,\lambda}(n,k).
	\end{displaymath}
\end{theorem}
From \eqref{6}, we get
\begin{align}
\sum_{n=k}^{\infty}S_{2,\lambda}(n,k)\frac{t^{n}}{n!}\ &=\ \frac{1}{k!}\big(e_{\lambda}(t)-1\big)^{k}\ =\ \frac{1}{k!}\sum_{l=0}^{k}\binom{k}{l}(-1)^{k-l}e_{\lambda}^{l}(t) \label{25} \\
&=\ \sum_{n=0}^{\infty}\bigg(\frac{1}{k!}\sum_{l=0}^{k}\binom{k}{l}(-1)^{k-l}(l)_{n,\lambda}\bigg)\frac{t^{n}}{n!}. \nonumber
\end{align}
By \eqref{25}, we get
\begin{equation}
\sum_{l=0}^{k}\binom{k}{l}(-1)^{k-l}(l)_{n,\lambda}=\left\{\begin{array}{cc}
	k!S_{2,\lambda}(n,k), & \textrm{if $n\ge k$}, \\
	0, & \textrm{otherwise.}
\end{array}\right.	\label{26}
\end{equation}
Let $\triangle$ be a difference operator with $\triangle f(x)=f(x+1)-f(x)$. Then we have
\begin{displaymath}
	\triangle^{n}f(x)=\sum_{k=0}^{n}\binom{n}{k}(-1)^{n-k}f(x+k).
\end{displaymath}
From \eqref{26}, we have
\begin{equation}
k!S_{2,\lambda}(n,k)=\triangle^{k}(0)_{n,\lambda},\quad(n,k\ge 0, n\ge k).\label{27}
\end{equation}
In light of \eqref{12}, we may consider the {\it{degenerate type Eulerian numbers}} given by
\begin{equation}
(-1)^{n-m}\eulerian{n}{m}_{\lambda}=\sum_{k=0}^{n-m}\lambda^{k}(1)_{k+1,1/\lambda}\binom{n-k}{m}\frac{\triangle^{k}(0)_{n,\lambda}}{k!}.\label{28}
\end{equation}
By \eqref{27} and \eqref{28}, we get
\begin{equation}
(-1)^{n-m}\eulerian{n}{m}_{\lambda}=\sum_{k=0}^{n-m}\lambda^{k}(1)_{k+1,1/\lambda}\binom{n-k}{m}S_{2,\lambda}(n,k).\label{29}
\end{equation}
We observe that
\begin{align}
\sum_{k=0}^{n}\lambda^{k}(1)_{k+1,1/\lambda}S_{2,\lambda}(n,k)(t+1)^{n-k}\ &=\ \sum_{k=0}^{n}\lambda^{k}(1)_{k+1,1/\lambda}S_{2,\lambda}(n,k)\sum_{m=0}^{n-k}\binom{n-k}{m}t^{m}\label{30} \\
&=\ \sum_{m=0}^{n}\bigg(\sum_{k=0}^{n-m}\lambda^{k}(1)_{k+1,1/\lambda}S_{2,\lambda}(n,k)\binom{n-k}{m}\bigg)t^{m} \nonumber \\
&=\ \sum_{m=0}^{n}(-1)^{n-m}\eulerian{n}{m}_{\lambda}t^{m}.	\nonumber
\end{align}
From \eqref{30} and Theorem 3, we note that
\begin{align}
\beta_{n,\lambda}^{(p)}\ &=\ \sum_{k=0}^{n}\lambda^{k}(1)_{k+1,1/\lambda}\binom{p+k+1}{k}^{-1}S_{2,\lambda}(n,k) \label{31} \\
&=\ (p+1)\sum_{k=0}^{n}\lambda^{k}(1)_{k+1,1/\lambda}S_{2,\lambda}(n,k)\int_{0}^{1}t^{p}(1-t)^{k}dt\nonumber\\
&=\ (p+1)\int_{0}^{1}\sum_{k=0}^{n}\lambda^{k}(1-t)^{n}t^{p}(1)_{k+1,1/\lambda}S_{2,\lambda}(n,k)\bigg(1+\frac{t}{1-t}\bigg)^{n-k}dt \nonumber \\
&=\ (p+1)\int_{0}^{1}(1-t)^{n}t^{p}\sum_{k=0}^{n}\eulerian{n}{k}_{\lambda}(-1)^{n-k}\bigg(\frac{t}{1-t}\bigg)^{k}dt \nonumber \\
&=\ (p+1)\sum_{k=0}^{n}\eulerian{n}{k}_{\lambda}(-1)^{n-k}\int_{0}^{1} (1-t)^{n-k}t^{p+k}dt\nonumber \\
&=\ (p+1)\sum_{k=0}^{n}\eulerian{n}{k}_{\lambda}(-1)^{n-k}\frac{(n-k)!(p+k)!}{(p+n+1)!}\nonumber \\
&=\ \frac{p+1}{n+p+1}\sum_{k=0}^{n}\eulerian{n}{k}_{\lambda}(-1)^{n-k}\binom{p+n}{p+k}^{-1}. \nonumber
\end{align}
Therefore, by \eqref{31}, we obtain the following theorem.
\begin{theorem}
	For $n,p\ge 0$, we have
	\begin{displaymath}
		\beta_{n,\lambda}^{(p)}= \frac{p+1}{n+p+1}\sum_{k=0}^{n}\eulerian{n}{k}_{\lambda}(-1)^{n-k}\binom{p+n}{p+k}^{-1}.
	\end{displaymath}
\end{theorem}
Let $r$ be a positive integer. The unsigned $r$-Stirling number of the first kind ${n \brack k}_{r}$ is the number of permutations of the set $[n]=\{1,2,3,\dots,n\}$ with exactly $k$ disjoint cycles in such a way that the numbers $1,2,3,\dots,r$ are in distinct cycles, while the $r$-Stirling number of the second kind ${n \brace k}_{r}$ counts the number of partitions of the set $[n]$ into $k$ non-empty disjoint subsets in such a way that the numbers $1,2,3,\dots,r$ are in distinct subsets. In [13], Kim-Kim-Lee-Park introduced the unsigned degenerate $r$-Stirling numbers of the first kind ${n \brack k}_{r,\lambda}$ as a degenerate version of ${n \brack k}_{r}$ and the degenerate $r$-Stirling number of the second kind ${n \brace k}_{r,\lambda}$ as a degenerate version of ${n \brace k}_{r}$. It is known that the degenerate $r$-Stirling numbers of the second kind are given by
\begin{equation}
(x+r)_{n,\lambda}=\sum_{k=0}^{n}{n+r \brace k+r}_{r,\lambda}(x)_{k},\quad(n \ge 1). \label{32}
\end{equation}
From \eqref{32}, we note that
\begin{equation}
\frac{1}{k!}\big(e_{\lambda}(t)-1\big)^{k}e_{\lambda}^{r}(t)=\sum_{n=k}^{\infty}{n+r \brace k+r}_{r,\lambda}\frac{t^{n}}{n!},\quad(k\ge 0,\ r \ge 1). \label{33}
\end{equation}
By the Euler's transformation formula in \eqref{9}, we get
\begin{align}
\sum_{n=0}^{\infty}\beta_{n,\lambda}^{(p)}\frac{t^{n}}{n!}\ &=\ \pFq{2}{1}{1-\lambda,1}{p+2}{1-e_{\lambda}(t)} \label{34} \\
&=\ e_{\lambda}^{p+\lambda}(t)\sum_{k=0}^{\infty}\frac{\langle p+1+\lambda\rangle_{\lambda}\langle p+1\rangle_{k}}{\langle p+2\rangle_{\lambda}}\frac{(1-e_{\lambda}(t))^{k}}{k!}\nonumber \\
&=\ \frac{p+1}{\lambda^{p}\langle 1\rangle_{p+1,1/\lambda}}\sum_{k=0}^{\infty}\frac{\lambda^{k+p}\langle 1\rangle_{p+k+1,1/\lambda}}{p+k+1}\frac{(-1)^{k}}{k!}\big(e_{\lambda}(t)-1\big)^{k}e_{\lambda}^{p+k}(t)\nonumber\\
&=\ \frac{p+1}{\lambda^{p}\langle 1\rangle_{p+1,1/\lambda}}\sum_{k=0}^{\infty}\frac{\lambda^{k+p}\langle 1\rangle_{p+k+1,1/\lambda}}{p+k+1}(-1)^{k}\sum_{m=k}^{\infty}{m+p \brace k+p}_{p,\lambda}\frac{t^{m}}{m!}\sum_{l=0}^{\infty}(k)_{l,\lambda}\frac{t^l}{l!} \nonumber  \\
&=\ \sum_{m=0}^{\infty}\frac{p+1}{\lambda^{p}\langle 1\rangle_{p+1,1/\lambda}}\sum_{k=0}^{m}\frac{\lambda^{k+p}\langle 1\rangle_{p+k+1,1/\lambda}}{p+k+1}(-1)^{k}{m+p \brace k+p}_{p,\lambda}\frac{t^{m}}{m!}\sum_{l=0}^{\infty}(k)_{l,\lambda}\frac{t^l}{l!}\nonumber \\
&=\ \sum_{n=0}^{\infty}\sum_{m=0}^{n}\binom{n}{m}\frac{p+1}{\langle 1\rangle_{p+1,1/\lambda}}\sum_{k=0}^{m}\frac{(-\lambda)^{k}}{p+k+1} \langle 1\rangle_{p+k+1,1/\lambda}{m+p \brace k+p}_{p,\lambda}(k)_{n-m,\lambda}\frac{t^{n}}{n!}, \nonumber
\end{align}
where $\langle x\rangle_{0,\lambda}=1,\ \langle x\rangle_{n,\lambda}=x(x+\lambda)\cdots(x+(n-1)\lambda),\ (n\ge 1)$. \par

\vspace{0.1cm}

Therefore, we obtain the following theorem.
\begin{theorem}
For $n\ge 1$ and $p\ge 0$, we have
\begin{equation*}
\beta_{n,\lambda}^{(p)}=\frac{p+1}{\langle 1\rangle_{p+1,1/\lambda}} \sum_{m=0}^{n}\sum_{k=0}^{m}\binom{n}{m}\frac{(-\lambda)^{k}}{p+k+1} \langle 1\rangle_{p+k+1,1/\lambda}{m+p \brace k+p}_{p,\lambda}(k)_{n-m,\lambda}.
\end{equation*}
\end{theorem}
Note that
\begin{displaymath}
\lim_{\lambda\rightarrow 0}\beta_{n,\lambda}^{(p)}=\frac{p+1}{p!}\sum_{m=0}^{n}\sum_{k=0}^{m}\binom{n}{m}(-1)^{k}\frac{(p+k)!}{p+k+1}{m+p \brace k+p}_{p}(k)_{n-m}.
\end{displaymath}
From Theorem 3, we have
\begin{align}
\sum_{n=0}^{\infty}\beta_{n,\lambda}^{(p)}\frac{t^{n}}{n!}\ &=\ \sum_{n=0}^{\infty}\bigg(\sum_{k=0}^{n}\frac{\lambda^{k}(1)_{k+1,1/\lambda}}{\binom{p+k+1}{p+1}}S_{2,\lambda}(n,k)\bigg)\frac{t^{n}}{n!}\label{35} \\
&=\ \sum_{k=0}^{\infty}\frac{\lambda^{k}(1)_{k+1,1/\lambda}}{\binom{p+k+1}{p+1}}\frac{1}{k!}\big(e_{\lambda}(t)-1\big)^{k}\nonumber \\
&=\ (p+1)\sum_{k=0}^{\infty}\frac{p!k!}{(k+p+1)!}\lambda^{k}(1)_{k+1,1/\lambda}\frac{1}{k!}\big(e_{\lambda}(t)-1\big)^{k}	\nonumber \\
&=\ (p+1)\sum_{k=0}^{\infty}\frac{(-1)^{k}\lambda^{k}(1)_{k+1,1/\lambda}}{k!}\big(1-e_{\lambda}(t)\big)^{k}\int_{0}^{1}(1-x)^{p}x^{k}dx\nonumber \\
&=\ (p+1)\sum_{k=0}^{\infty}(-1)^{k}\binom{\lambda-1}{k}\big(1-e_{\lambda}(t)\big)^{k}\int_{0}^{1}(1-x)^{p}x^{k}dx\nonumber \\
&=\ (p+1)\int_{0}^{1}(1-x)^{p}\big(1-x(1-e_{\lambda}(t))\big)^{\lambda-1}dx. \nonumber
\end{align}
Therefore, we obtain the following theorem.
\begin{theorem}
	For $p\ge 0$, we have
	\begin{displaymath}
		\sum_{n=0}^{\infty}\beta_{n,\lambda}^{(p)}\frac{t^{n}}{n!}=(p+1)\int_{0}^{1}(1-x)^{p}\big(1-x(1-e_{\lambda}(t))\big)^{\lambda-1}dx.
	\end{displaymath}
\end{theorem}

\section{Generalized degenerate Bernoulli polynomials}
In this section, we consider the {\it{generalized degenerate Bernoulli polynomials}} which are derived from the Gauss hypergeometric function.
In light of \eqref{21}, we define the generalized degenerate Bernoulli polynomials by
\begin{equation}
\sum_{n=0}^{\infty}\beta_{n,\lambda}^{(p)}(x)\frac{t^{n}}{n!}=\pFq{2}{1}{1-\lambda,1}{p+2}{1-e_{\lambda}(t)}e_{\lambda}^{x}(t). \label{36}
\end{equation}
When $x=0$, $\beta_{n,\lambda}^{(p)}(0)=\beta_{n,\lambda}^{(p)}$, $(n\ge 0)$. Thus, by \eqref{36}, we get
\begin{align}
\sum_{n=0}^{\infty}\beta_{n,\lambda}^{(p)}(x)\frac{t^{n}}{n!}\ &=\ \pFq{2}{1}{1-\lambda,1}{p+2}{1-e_{\lambda}(t)}e_{\lambda}^{x}(t) \label{37} \\
&=\ \sum_{l=0}^{\infty}\beta_{l,\lambda}^{(p)}\frac{t^{l}}{l!}\sum_{m=0}^{\infty}(x)_{m,\lambda}\frac{t^{m}}{m!}\nonumber \\
&=\ \sum_{n=0}^{\infty}\bigg(\sum_{l=0}^{n}\binom{n}{l}\beta_{l,\lambda}^{(p)}(x)_{n-l,\lambda}\bigg)\frac{t^{n}}{n!}. \nonumber
\end{align}
Therefore, by comparing the coefficients on both sides of \eqref{37}, we obtain the following theorem.
\begin{theorem}
For $n\ge 0$, we have
\begin{displaymath}
\beta_{n,\lambda}^{(p)}(x)=\sum_{l=0}^{n}\binom{n}{l}\beta_{l,\lambda}^{(p)}(x)_{n-l,\lambda}.
\end{displaymath}
\end{theorem}
From \eqref{36}, we note that
\begin{align*}
	\sum_{n=1}^{\infty}\frac{d}{dx}\beta_{n,\lambda}^{(p)}(x)\frac{t^{n}}{n!}\ &=\ \pFq{2}{1}{1-\lambda,1}{p+2}{1-e_{\lambda}(t)}\frac{d}{dx}e_{\lambda}^{x}(t) \\
	&=\ \frac{1}{\lambda}\log(1+\lambda t)\pFq{2}{1}{1-\lambda,1}{p+2}{1-e_{\lambda}(t)}e_{\lambda}^{x}(t) \\
	&=\ \frac{1}{\lambda}\sum_{l=1}^{\infty}\frac{(-1)^{l-1}\lambda^{l}}{l}t^{l}\sum_{m=0}^{\infty}\beta_{m,\lambda}^{(p)}(x)\frac{t^{m}}{m!} \\
	&=\ \sum_{n=1}^{\infty}\bigg(\sum_{l=1}^{n}\frac{(-\lambda)^{l-1}}{l}\frac{n!\beta_{n-l,\lambda}^{(p)}(x)}{(n-l)!}\bigg)\frac{t^{n}}{n!}
\end{align*}
Thus, we have
\begin{displaymath}
	\frac{d}{dx}\beta_{n,\lambda}^{(p)}(x)= \sum_{l=1}^{n}(-\lambda)^{l-1}(l-1)!\binom{n}{l}\beta_{n-l,\lambda}^{(p)}(x).
\end{displaymath}
\begin{proposition}
	For $n\ge 1$, we have
	\begin{displaymath}
		\frac{d}{dx}\beta_{n,\lambda}^{(p)}(x)= \sum_{l=1}^{n}(-\lambda)^{l-1}(l-1)!\binom{n}{l}\beta_{n-l,\lambda}^{(p)}(x).
	\end{displaymath}
\end{proposition}
By \eqref{14}, we easily get
\begin{align*}
	\sum_{n=k}^{\infty}S_{2,\lambda}(n,k|x)\frac{t^{n}}{n!}\ &=\ \frac{1}{k!}\big(e_{\lambda}(t)-1\big)^{k}e_{\lambda}^{x}(t) \\
	&=\ \frac{1}{k!}\sum_{l=0}^{k}\binom{k}{l}(-1)^{k-l}e_{\lambda}^{l+x}(t)  \\
	&=\ \sum_{n=0}^{\infty}\bigg(\frac{1}{k!}\sum_{l=0}^{k}(-1)^{k-l}(l+x)_{n,\lambda}\bigg)\frac{t^{n}}{n!}.
\end{align*}
Thus we have
\begin{equation}
\frac{1}{k!}\sum_{l=0}^{k}\binom{k}{l}(-1)^{k-l}(l+x)_{n,\lambda}=\left\{\begin{array}{ccc}
	S_{2,\lambda}(n,k|x), & \textrm{if $n\ge k$,}\\
	0, & \textrm{otherwise.}
\end{array}\right.	\label{38}
\end{equation}
From \eqref{38}, we note that
\begin{displaymath}
	S_{2,\lambda}(n,k|x)=\frac{1}{k!}\triangle^{k}(x)_{n,\lambda},\quad(n\ge k).
\end{displaymath}
\begin{lemma}
	For $n,k\ge 0$ with $n\ge k$, we have
	\begin{displaymath}
	S_{2,\lambda}(n,k|x)=\frac{1}{k!}\triangle^{k}(x)_{n,\lambda},\quad(n\ge k).
\end{displaymath}
\end{lemma}
Now, we observe that
\begin{align}
\sum_{n=0}^{\infty}\beta_{n,\lambda}^{(p)}(x)\frac{t^{n}}{n!}\ &=\ \sum_{k=0}^{\infty}\frac{(p+1)!k!}{(p+k+1)!}\lambda^{k}(1)_{k+1,1/\lambda}\frac{1}{k!}\big(e_{\lambda}(t)-1\big)^{k}e_{\lambda}^{x}(t) \label{39} \\
&=\ \sum_{k=0}^{\infty}\frac{(1)_{k+1,1/\lambda}\lambda^{k}}{\binom{p+k+1}{p+1}}\sum_{n=k}^{\infty}S_{2,\lambda}(n,k|x)\frac{t^{n}}{n!} \nonumber \\
&=\ \sum_{n=0}^{\infty}\bigg(\sum_{k=0}^{n}\frac{(1)_{k+1,1/\lambda}\lambda^{k}}{\binom{p+k+1}{p+1}}S_{2,\lambda}(n,k|x)\bigg)\frac{t^{n}}{n!}.\nonumber
\end{align}
Therefore, by \eqref{39}, we obtain the following theorem.
\begin{theorem}
For $n\ge 0$, we have
\begin{displaymath}
\beta_{n,\lambda}^{(p)}(x)= \sum_{k=0}^{n}\frac{(1)_{k+1,1/\lambda}\lambda^{k}}{\binom{p+k+1}{p+1}}S_{2,\lambda}(n,k|x).
\end{displaymath}
\end{theorem}
\begin{remark}
Let $p$ be a nonnegative integer. Then, by Theorem 7 and \eqref{36}, we easily get
\begin{align*}
&\ \beta_{n,\lambda}^{(p)}(x+y)=\sum_{k=0}^{n}\binom{n}{k}\beta_{k,\lambda}^{(p)}(x)(y)_{n-k,\lambda},\quad(n \ge 0),\nonumber \\
&\ \beta_{n,\lambda}^{(p)}(x+1)-\beta_{n,\lambda}^{(p)}(x)\ =\ \sum_{k=0}^{n-1}\binom{n}{k}\beta_{k,\lambda}^{(p)}(x)(1)_{n-k,\lambda},\quad(n \ge 1),\nonumber \\
&\ \beta_{n,\lambda}^{(p)}(mx)\ =\ \sum_{k=0}^{n}\binom{n}{k}\beta_{k,\lambda}^{(p)}(x)(m-1)^{n-k}(x)_{n-k,\lambda/m-1},\quad(n \ge 0, m \ge 2). \nonumber
\end{align*}
\end{remark}

\section{Conclusion}
In recent years, degenerate versions of many special polynomials and numbers have been investigated by means of various different tools including generating functions, combinatorial methods, umbral calculus, $p$-adic analysis, differential equations, special functions, probability theory and analytic number theory. Studying degenerate versions of some special polynomials and numbers, which was initiated by Carlitz in [3], has yielded many interesting arithmetic and combinatorial results (see [7-13 and references therein]) and has potential to find many applications in diverse areas.\par
A new family of $p$-Bernoulli numbers and polynomials, which reduce to the classical Bernoulli numbers and polynomials for $p=0$, was introduced by Rahmani in [14] by means of the Gauss hypergeometric function. Motivated by that paper, we were interested in finding a degenerate version of those numbers and polynomials. Indeed, the generalized degenerate Bernoulli numbers and polynomials, which reduce to the Carlitz degenerate Bernoulli numbers and polynomials for $p=0$, were introduced again in terms of the Gauss hypergeometric function. In addition, the degenerate type Eulerian numbers was introduced as a degenerate version of Eulerian numbers. \par
In this paper, we expressed the generalized degenerate Bernoulli numbers in terms of the degenerate Stirling numbers of the second kind, of the degenerate type Eulerian numbers, of the degenerate $p$-Stirling numbers of the second kind and of an integral on the unit interval. In dddition, we represented the generalized degenerate Bernoulli polynomials in terms of the degenerate Stirling polynomials of the second kind. \par
It is one of our future projects to continue pursuing this line of research. Namely, by studying degenerate versions of some special polynomials and numbers, we want to find their applications in mathematics, science and engineering.

\vspace{0.2in}

{\bf Acknowledgements}

The authors would like to thank the reviewers for their valuable comments and suggestions and Jangjeon Institute for Mathematical Science for the support of this research.

\vspace{0.2in}

{\bf Funding}

Not applicable.

\vspace{0.2in}

{\bf Availability of data and materials}

Not applicable.

\vspace{0.2in}

{\bf Competing interests}

The authors declare that they have no conflicts of interest.

\vspace{0.2in}

{\bf Authors’ contributions}

TK and DSK conceived of the framework and structured the whole paper; DSK and TK wrote the paper; HL typed;  LCJ, and HYK
checked the results of the paper; DSK and TK completed the revision of the paper. All authors have read and approved
the final version of the manuscript.

\end{document}